\documentclass[12pt]{article}
\pdfoutput=1

\usepackage{a4wide}
\usepackage{amsmath}
\usepackage{amssymb}
\usepackage{amsthm}
\usepackage{picins}
\usepackage{graphicx}
\usepackage{wrapfig}
\usepackage{color}
\usepackage[small,bf]{caption}
\usepackage{ifpdf}
\ifpdf
\fi

\newtheorem{theorem}{Theorem}

\newtheorem{lemma}{Lemma}
\newtheorem{corollary}{Corollary}
\newtheorem{definition}{Definition}
\newtheorem{remark}{Remark}

\theoremstyle{definition}
\newtheorem{example}{Example}

\newcommand{\symm}{\bigtriangleup}

\begin{document}

\begin{center}
{\Huge On the difference between solutions of discrete tomography problems}

\bigskip

{\Large Birgit van Dalen}

\textit{Mathematisch Instituut, Universiteit Leiden, Niels Bohrweg 1, 2333 CA Leiden, The Netherlands \\ dalen@math.leidenuniv.nl}

\end{center}

{\small \textbf{Abstract:} We consider the problem of reconstructing binary images from their horizontal and vertical projections. We present a condition that the projections must necessarily satisfy when there exist two disjoint reconstructions from those projections. More generally, we derive an upper bound on the symmetric difference of two reconstructions from the same projections. We also consider two reconstructions from two different sets of projections and prove an upper bound on the symmetric difference in this case.}

\section{Introduction}\label{introduction}

Discrete tomography is concerned with problems such as reconstructing binary images on a lattice from given projections in lattice directions \cite{boek}. Each point of a binary image has a value equal to zero or one. The line sum of a line through the image is the sum of the values of the points on this line. The projection of the image in a certain lattice direction consists of all the line sums of the lines through the image in this direction.

Discrete tomography has applications in a wide range of fields. The most important are electron microscopy \cite{gold} and medical imaging \cite{medicaldt, medical}, but there are also applications in nuclear science \cite{radio2, radio} and various other fields \cite{propellor, building}.

For any set of directions, it is possible to construct binary images that are not uniquely determined by their projections in those directions \cite[Theorem 4.3.1]{boek}. The problem of deciding whether an image is uniquely determined by its projections and the problem of reconstructing it are NP-complete for any set of more than two directions \cite{gardner}. For exactly two directions, these problems can be solved in polynomial time. Already in 1957, Ryser described an algorithm to reconstruct binary images from their horizontal and vertical projections and characterised the set of projections that correspond to a unique binary image \cite{ryser}.

When a binary image is not uniquely determined by its projections, the reconstruction may not be equal to the original image. In practical applications, it is important to know whether the reconstruction is a good approximation of the original image. In some applications it may be sufficient to find a reconstruction of which a large part is guaranteed to belong to the original image. It is also interesting from a theoretical point of view to find bounds on how much two images with the same projections can differ, and to have conditions under which the two images can be completely disjoint.

There is a simple bound in the case of two directions: if the image is contained in an $m \times n$-rectangle and a certain row sum is equal to $a \geq \frac{1}{2}n$, then the difference in that row can be at most $2a-n$. If on the other hand a row sum is equal to $b < \frac{1}{2}n$, then the difference in the row can be at most $2b$. Summing over all $m$ rows gives an upper bound on the size of the symmetric difference of two different reconstructions. While this bound may be quite good in some special cases, it is not very good in general.

In this paper we will use a different approach, based on the work in \cite{birgit}. There the concept of staircases, introduced by Alpers \cite{alpersthesis}, was used to compare an arbitrary image to a uniquely determined image. Here we generalise this methode in order to compare two arbitrary binary images. We use a uniquely determined image that is as close as possible to the original image. We characterise such images in Theorem \ref{closestuniquesets}. We then consider two reconstructions of the same original image and prove bounds on the intersection and symmetric difference of the two reconstructions in Theorems \ref{dezelfde} and \ref{beta}. As a consequence of these results, we find a condition on the projections that must hold when the reconstruction and the original image are disjoint.

In Theorem \ref{metF1} we show that we can construct a uniquely determined image that is guaranteed to have a large intersection with the original image. As a consequence of these results, we find a condition on the projections that must hold when the reconstruction and the original image are disjoint. To complement this result, we state conditions under which no individual point must necessarily belong to the original image (these conditions are a direct consequence of a theorem by Anstee \cite{anstee}).
Finally, we will consider two different sets of horizontal and vertical projections and prove an upper bound for the difference between two reconstructions in this case.

\section{Notation}\label{notation}

Let $F$ be a finite subset of $\mathbb{Z}^2$ with characteristic function $\chi$. (That is, $\chi(x,y) = 1$ if and only if $(x,y) \in F$.) For $i \in \mathbb{Z}$, we define \emph{row} $i$ as the set $\{(x,y) \in \mathbb{Z}^2: x = i\}$. We call $i$ the index of the row. For $j \in \mathbb{Z}$, we define \emph{column} $j$ as the set $\{(x,y) \in \mathbb{Z}^2: y = j\}$. We call $j$ the index of the column. We order the rows and columns as is usual for a matrix: the row numbers increase when going downwards; the column numbers increase when going to the right.

The \emph{row sum} $r_i$ is the number of elements of $F$ in row $i$, that is $r_i = \sum_{j \in \mathbb{Z}} \chi(i,j)$. The \emph{column sum} $c_j$ of $F$ is the number of elements of $F$ in column $j$, that is $c_j = \sum_{i \in \mathbb{Z}} \chi(i,j)$. We refer to both row and column sums as the \emph{line sums} of $F$. We will usually only consider finite sequences $\{r_1, r_2, \ldots, r_m\}$ and $\{c_1, c_2, \ldots, c_n\}$ of row and column sums that contain all the nonzero line sums.

We call $F$ \emph{uniquely determined by its line sums} or simply \emph{uniquely determined} if the following property holds: if $F'$ is a subset of $\mathbb{Z}^2$ with exactly the same row and column sums as $F$, then $F' = F$. Suppose $F$ is uniquely determined and has row sums $r_1$, $r_2$, \ldots, $r_m$. For each $j$ with $1 \leq j \leq \max_i r_i$ we can count the number $\#\{l: r_l \geq j\}$ of row sums that are at least $j$. These numbers are exactly the non-zero column sums of $F$ (in some order). This is an immediate consequence of Ryser's theorem (\cite{ryser}, see also \cite[Theorem 1.7]{boek}).

Suppose we have two finite subsets $F_1$ and $F_2$ of $\mathbb{Z}^2$. For $h = 1, 2$ we denote the row and column sums of $F_h$ by $r_i^{(h)}$, $i \in \mathbb{Z}$, and $c_j^{(h)}$, $j \in \mathbb{Z}$, respectively. Define
\[
\alpha(F_1, F_2) = \frac{1}{2} \left( \sum_{j \in \mathbb{Z}} |c_j^{(1)} - c_j^{(2)}| + \sum_{i
\in \mathbb{Z}} |r_i^{(1)} - r_i^{(2)}| \right).
\]
Note that $\alpha(F_1, F_2)$ is always an integer, since $2 \alpha(F_1, F_2)$ is congruent to
\[
\sum_{j \in \mathbb{Z}} \left(c_j^{(1)} + c_j^{(2)}\right) + \sum_{i \in \mathbb{Z}} \left( r_i^{(1)} + r_i^{(2)} \right) = 2|F_1| + 2|F_2| \equiv 0 \mod 2.
\]
We will sometimes refer to $\sum_{j \in \mathbb{Z}} |c_j^{(1)} - c_j^{(2)}|$ as the \emph{difference in the column sums} and to $\sum_{i \in \mathbb{Z}} |r_i^{(1)} - r_i^{(2)}|$ as the \emph{difference in the row sums}.

In order to describe the symmetric difference between two sets $F_1$ and $F_2$, we use the notion of a staircase, first introduced by Alpers \cite{alpersthesis}.

\begin{definition}
A set of points $(p_1, p_2, \ldots, p_n)$ in $\mathbb{Z}^2$ is called a
\emph{staircase} if the following two conditions are satisfied:
\begin{itemize}
\item for each $i$ with $1 \leq i \leq n-1$ one of the points
$p_i$ and $p_{i+1}$ is an element of $F_1 \backslash F_2$ and the
other is an element of $F_2 \backslash F_1$;
\item either for all $i$ the points $p_{2i}$ and $p_{2i+1}$ are in
the same column and the points $p_{2i+1}$ and $p_{2i+2}$ are in
the same row, or for all $i$ the points $p_{2i}$ and $p_{2i+1}$
are in the same row and the points $p_{2i+1}$ and $p_{2i+2}$ are
in the same column.
\end{itemize}
\end{definition}

\section{Some lemmas}\label{lemmas}

We prove some lemmas that we will use later for our main results.

\begin{lemma}\label{rijtjes}
Let $a_1 \geq a_2 \geq \ldots \geq a_n$ be non-negative integers. Let $m \geq \max_j a_j$. For $1 \leq i \leq m$, define $b_i
= \# \{ j : a_j \geq i \}$. Then for $1 \leq j \leq n$ we have $a_j = \# \{ i : b_i \geq j\}$.
\end{lemma}

\begin{proof} We have $b_1 \geq b_2 \geq \ldots \geq b_m$. Hence for $1 \leq j \leq n$ we have
\[
\# \{ i : b_i \geq j\} = \max\{i : b_i \geq j\} = \max\{i : \max\{l: a_l \geq i\} \geq j\}.
\]
For a fixed $i$ we have
\[
\max\{l: a_l \geq i\} \geq j
\quad \Longleftrightarrow \quad a_j \geq i,
\]
hence
\[
\max\{i : \max\{l: a_l \geq i\} \geq j\} = \max\{ i : a_j \geq i \} = a_j.
\]
This completes the proof.
\end{proof}

\begin{lemma}\label{lemmauniek}
Let $F$ be a uniquely determined finite subset of $\mathbb{Z}^2$ with row sums $r_i$, $i \in \mathbb{Z}$, and column sums $c_j$, $j \in \mathbb{Z}$, respectively. If for integers $i_1$, $i_2$ and $j_0$ we have $(i_1,j_0) \in F$ and $(i_2,j_0) \not\in F$, then $r_{i_1} > r_{i_2}$.
\end{lemma}

\begin{proof}
As $F$ is uniquely determined, we have the following characterisation of its elements \cite[p. 17]{boek}: a point $(x,y)$ is an element of $F$ if and only if $r_x \geq \#\{l: c_l \geq c_y\}$. So if $(i_1,j_0) \in F$ and $(i_2,j_0) \not\in F$, we have $r_{i_1} \geq \#\{l: c_l \geq c_{j_0}\} > r_{i_2}$.
\end{proof}

Let $F_1$ and $F_2$ be finite subsets of $\mathbb{Z}^2$, such that $F_1$ is uniquely determined and $|F_1| = |F_2|$. Denote the row sums of $F_1$ by $r_i$, $i \in \mathbb{Z}$. Let $\alpha = \alpha(F_1, F_2)$. The symmetric difference $F_1 \symm F_2$ is the disjoint union of $\alpha$ staircases \cite{birgit}. Consider such a staircase with points $(x_1, y_1), (x_2, y_2), \ldots, (x_t,y_t) \in F_1 \backslash F_2$ and $(x_2, y_1), (x_3, y_2) \ldots, (x_t, y_{t-1}) \in F_2 \backslash F_1$. (The staircase may contain another point of $F_2 \backslash F_1$ in row $x_1$ and another one in column $y_t$, but this is irrelevant here.) By Lemma \ref{lemmauniek} we have
\[
r_{x_1} > r_{x_2} > \ldots > r_{x_t}.
\]
Hence the rows $x_1$, $x_2$, \ldots, $x_t$ of $F_1$ have pairwise different line sums.

\begin{lemma}\label{wortelgrens}
We have
\[
|F_1 \symm F_2| \leq \alpha \sqrt{8|F_1| + 1} - \alpha.
\]
\end{lemma}

\begin{proof}
Let $n$ be the largest positive integer such that $|F_1| \geq \frac{1}{2}n(n+1)$. Suppose $F_1$ has at least $n+1$ distinct positive row sums. Then
\[
|F_1| \geq 1 + 2 + \cdots + n + (n+1) = \frac{1}{2}(n+1)(n+2),
\]
which contradicts the maximality of $n$. So $F_1$ has at most $n$ distinct positive row sums. Any staircase of $F_1 \symm F_2$ therefore contains elements of $F_1 \backslash F_2$ in at most $n$ different rows. So the total number of elements of $F_1 \backslash F_2$ cannot exceed $\alpha n$. Hence $|F_1 \symm F_2| \leq 2 \alpha n$. On the other hand, we have $2 |F_1| \geq n^2+n = (n+\frac{1}{2})^2 - \frac{1}{4}$, thus $n \leq \sqrt{2|F_1| + \frac{1}{4}} - \frac{1}{2}$. We conclude
\[
|F_1 \symm F_2| \leq \alpha \sqrt{8|F_1| + 1} - \alpha.
\]
\end{proof}

\begin{remark}\label{wortelgrens2}
We will also use the slightly weaker estimate
\[
|F_1 \symm F_2| \leq 2 \alpha \sqrt{2|F_1|}.
\]
\end{remark}

\section{Uniquely determined neighbours}\label{neighbour}

Consider a set $F_2$ that is not uniquely determined by its line sums. We are interested in how close -- in some sense -- this set is to being uniquely determined. We define the distance between $F_2$ and a uniquely determined set $F_1$ as $\alpha(F_2, F_1)$. The smallest possible value of $\alpha(F_2, F_1)$ then indicates how close $F_2$ is to being uniquely determined. It turns out that we can characterise in a very simple way the sets $F_1$ for which $\alpha(F_2, F_1)$ is minimal.

\begin{theorem}\label{closestuniquesets}
Let $F_2$ be a finite subset of $\mathbb{Z}^2$ with non-zero row sums $r_1 \geq r_2 \geq \ldots \geq r_m$ and non-zero column sums $c_1 \geq c_2 \geq \ldots \geq c_n$. Put $a_j = \#\{i : r_i \geq j \}$, $1 \leq j \leq n$, and $b_i = \#\{j : c_j \geq i\}$, $1 \leq i \leq m$. Define $\alpha_0 = \min\{\alpha(F_2, F) : F \text{ is a uniquely determined set}\}$. Let $F_1$ be a uniquely determined set with row sums $u_1 \geq u_2 \geq \ldots$, and column sums $v_1 \geq v_2 \geq \ldots$. Then the following conditions are equivalent:
\begin{itemize}
\item[(i)] $\alpha(F_2, F_1) = \alpha_0$,
\item[(ii)] for all $j \geq 1$ we have $\begin{cases}\min(a_j, c_j) \leq v_j \leq \max(a_j, c_j) & \text{if $1 \leq j \leq n$,} \\ v_j = 0 & \text{otherwise,} \end{cases}$
\item[(iii)] for all $i \geq 1$ we have $\begin{cases}\min(b_i, r_i) \leq u_i \leq \max(b_i, r_i) & \text{if $1 \leq i \leq m$,} \\ u_i = 0 & \text{otherwise.} \end{cases}$
\end{itemize}
\end{theorem}

\begin{proof}We will prove the equivalence of $(i)$ and $(ii)$. By symmetry the equivalence of $(i)$ and $(iii)$ then follows. During the proof, we will use several times the fact that $u_i = \#\{j : v_j \geq i\}$, $i \geq 1$, as $F_1$ is uniquely determined (see Section \ref{notation}).

$(i) \Rightarrow (ii)$. Suppose $F_1$ does not satisfy $(ii)$. Then either $v_j \neq 0$ for some $j > n$, or $v_j < \min(a_j, c_j)$ for some $j$ with $1 \leq j \leq n$, or $v_j > \max(a_j,c_j)$ for some $j$ with $1 \leq j \leq n$. In each of those three cases we will prove that there exists a uniquely determined set $F_1'$ such that $\alpha(F_2, F_1') < \alpha(F_2, F_1)$, which implies that $F_1$ does not satisfy $(i)$.

Case 1: there is an $l >n$ such that $v_l \neq 0$. As for all $j \in \{1, 2, \ldots, n\}$ we have $v_j \geq v_l$, we must have $u_{v_l} = \#\{j: v_j \geq v_l\} \geq n+1$. Now consider the set $F_1'$ with the same row and column sums as $F_1$, except that the column sum with index $l$ is exactly 1 smaller and the row sum with index $v_l$ is exactly 1 smaller. Note that this set is uniquely determined. Since either $v_l > m$ (so $r_{v_l}$ does not exist) or $r_{v_l} \leq n$, the difference in the row sums of $F_1'$ and $F_2$ is 1 less than the difference in the row sums of $F_1$ and $F_2$. The same holds for the differences in the column sums. So $\alpha(F_2, F_1') < \alpha(F_2, F_1)$.

Case 2: there is a $k \in \{1, 2, \ldots, n\}$ such that $v_k < \min(a_k, c_k)$. Assume that $k$ is the smallest positive integer with this property. Define $F_1'$ such that its row sums $u_i'$ and column sums $v_j'$ are as follows:
\[
u_i' = \begin{cases} u_i + 1 & \text{if $i = v_k+1$,} \\ u_i & \text{otherwise,} \end{cases}
\]
\[
v_j' = \begin{cases} v_k +1 & \text{if $j=k$,} \\ v_j & \text{otherwise.} \end{cases}
\]
If $k=1$, then the column sums of $F_1'$ are obviously non-increasing. If $k \geq 2$, then
\[
v_{k-1} \geq \min(a_{k-1}, c_{k-1}) \geq \min(a_k, c_k) > v_k,
\]
so $v_{k-1}' = v_{k-1} \geq v_k + 1 = v_k'$, hence the column sums are non-increasing in this case as well. For the row sums we have $u_i' = \#\{j : v_j' \geq i\}$, which shows that the row sums are non-increasing and that $F_1'$ is uniquely determined. \\ Clearly, the difference in the column sums has decreased by 1 when changing from $F_1$ to $F_1'$. The difference in the row sums has changed by $|u_{v_k+1} + 1 - r_{v_k+1}| - |u_{v_k+1} - r_{v_k+1}|$. We have $u_{v_k+1} = \#\{j : v_j \geq v_k + 1\} < k$. By Lemma \ref{rijtjes} we have $r_{v_k+1} = \#\{j : a_j \geq v_k + 1\}$ and therefore $r_{v_k+1} \geq k$, using $a_k \geq \min(a_k, c_k) \geq v_k + 1$. Hence $u_{v_k+1} < r_{v_k+1}$ and therefore the difference in the row sums has decreased by 1. So $\alpha(F_2, F_1') < \alpha(F_2, F_1)$.

Case 3: there is a $k \in \{1, 2, \ldots, n\}$ such that $v_k > \max(a_k, c_k)$. This is analogous to Case 2.

$(ii) \Rightarrow (i)$. Suppose $F_1$ satisfies $(ii)$. Consider the uniquely determined set with column sums $\min(a_j, c_j)$, $1 \leq j \leq n$, and non-increasing row sums. Then we can build $F_1$ starting from this set by adding new points one by one. Starting in the column with index 1, we add points to each column until there are $v_j$ points in column $j$. The points added in column $j$ are in rows $\min(a_j, c_j) + 1$, \ldots, $v_j$ in that order. In this way, in every step the constructed set has non-increasing row and column sums and is uniquely determined. We will prove that the value of $\alpha$ does not change in each step, which implies that the value of $\alpha$ of the set we started with is equal to $\alpha(F_2, F_1)$. That proves that all sets $F_1$ satisfying $(ii)$ have the same value $\alpha(F_2, F_1)$. This must then be the minimal value $\alpha_0$, since we proved in the first part that the minimal value occurs among the sets $F_1$ satisfying $(ii)$.

Now assume that $F_1$ satisfies $(ii)$ and let $k$ be such that $v_k < \max(a_k, c_k)$ and if $k \geq 2$, then $v_k < v_{k-1}$. It suffices to prove that if we add the point $(k, v_k+1)$ to $F_1$, then the value of $\alpha$ does not change. (Whenever we add a point in the procedure described above, the conditions $v_k < \max(a_k,c_k)$ and $v_k < v_{k-1}$ hold.) So define $F_1'$ as the uniquely determined set with row sums $u_i'$ and column sums $v_j'$ satisfying
\[
u_i' = \begin{cases} u_i + 1 & \text{if $i = v_k+1$,} \\ u_i & \text{otherwise,} \end{cases}
\]
\[
v_j' = \begin{cases} v_k +1 & \text{if $j=k$,} \\ v_j & \text{otherwise.} \end{cases}
\]
We will prove that $\alpha(F_2, F_1')= \alpha(F_2, F_1)$. We distinguish two cases.

Case 1: $a_k \leq v_k < c_k$. By changing from $F_1$ to $F_1'$ the difference in the column sums has decreased by 1. We have $u_{v_k+1} = \#\{j: v_j \geq v_k + 1\} = k-1$, as either $k=1$ or $v_{k-1} \geq v_k + 1$. Also, by Lemma \ref{rijtjes} we have $r_{v_k+1} = \#\{j: a_j \geq v_k + 1\} \leq k-1$, since $a_k < v_k + 1$. So $u_{v_k+1} \geq r_{v_k+1}$, which shows that the difference in the row sums has increased by 1. Hence $\alpha(F_2, F_1')= \alpha(F_2, F_1)$.

Case 2: $c_k \leq v_k < a_k$. By changing from $F_1$ to $F_1'$ the difference in the column sums has increased by 1. We have $u_{v_k+1} = k-1$ as in Case 1. Also, by Lemma \ref{rijtjes} we have $r_{v_k+1} = \#\{j: a_j \geq v_k + 1\} \geq k$, since $a_k \geq v_k+1$. So $u_{v_k+1} < r_{v_k+1}$, which shows that the difference in the row sums has decreased by 1. Hence $\alpha(F_2, F_1')= \alpha(F_2, F_1)$.

This completes the proof of the theorem.
\end{proof}

\begin{remark}
We can always permute the rows and columns such that the row and column sums of $F_2$ are non-increasing, so this condition in the above theorem is not a restriction. However, the monotony of the line sums of $F_1$ is a slight restriction. There may be a uniquely determined set $F_1$ satisfying $\alpha(F_2, F_1) = \alpha_0$ while its row and column sums are not non-increasing. However, reordering the row and column sums so that they are non-increasing never increases the differences with the row and column sums of $F_2$. So define in that case a set $F_1'$ with the same row and column sums as $F_1$, except that the line sums of $F_1'$ are ordered non-increasingly. Then $\alpha(F_2, F_1') = \alpha(F_2, F_1) = \alpha_0$, so $F_1'$ satisfies the conditions of the theorem and $(i)$ and therefore satisfies $(ii)$ and $(iii)$.
\end{remark}

Let $F_2$ be a set with row sums $r_1$, $r_2$, \ldots, $r_m$ and column sums $c_1$, $c_2$, \ldots, $c_n$, not necessarily non-increasing. Let $\sigma$ be a permutation of $\{1, 2, \ldots, n\}$ such that $c_{\sigma(1)} \geq c_{\sigma(2)} \geq \ldots \geq c_{\sigma(n)}$. Consider the uniquely determined set $F_1$ with row sums $u_1 = r_1$, $u_2 = r_2$, \ldots, $u_m = r_m$ and column sums $v_1$, $v_2$, \ldots, $v_n$ such that $v_{\sigma(1)} \geq v_{\sigma(2)} \geq \ldots \geq v_{\sigma(n)}$. According to Theorem \ref{closestuniquesets} we have $\alpha(F_2, F_1) = \alpha_0$, where $\alpha_0 = \min\{\alpha(F_2, F) : F \text{ is a uniquely determined set}\}$. Such a set $F_1$ we call a \emph{uniquely determined neighbour of $F_2$}. Note that $F_2$ may have more than one uniquely determined neighbour, as there may be more possibilities for $\sigma$ if some of the column sums of $F_2$ are equal. Also note that if $F_3$ is another set with row sums $r_1$, $r_2$, \ldots, $r_m$ and column sums $c_1$, $c_2$, \ldots, $c_n$, then $F_1$ is a uniquely determined neighbour of $F_3$ if and only if it is a uniquely determined neighbour of $F_2$.

It is easy to compute the line sums of a uniquely determined neighbour of $F_2$ and hence it is easy to find $\alpha_0$.

\section{Sets with equal line sums}\label{equallinesums}

Consider a set $F_2$ that is not uniquely determined by its line sums. When attempting to reconstruct $F_2$ from its line sums, one may end up with a different set $F_3$ that has the same line sums as $F_2$. It is interesting to know whether $F_3$ is a good approximation of $F_2$ or not. In some cases, $F_3$ may be disjoint from $F_2$, but in other cases, $F_2$ and $F_3$ must have a large intersection. We shall derive an upper bound on $F_2 \symm F_3$ that depends on the size of $F_2$ and on how close $F_2$ is to being uniquely determined, in the sense of the previous section. Both parameters can easily be computed from the line sums of $F_2$.

\begin{theorem}\label{dezelfde}
Let $F_2$ and $F_3$ be finite subsets of $\mathbb{Z}^2$ with the same line sums. Let $F_1$ be a uniquely determined neighbour of $F_2$ and $F_3$. Put $\alpha = \alpha(F_2, F_1)$. Then
\[
|F_2 \symm F_3| \leq 2 \alpha \sqrt{8|F_2| + 1} - 2 \alpha.
\]
\end{theorem}

\begin{proof}
By Lemma \ref{wortelgrens} we have $\alpha \sqrt{8|F_2| + 1} - \alpha$ as an upper bound for both $|F_1 \symm F_2|$ and $|F_1 \symm F_3|$. Hence
\[
|F_2 \symm F_3| \leq |F_1 \symm F_2| + |F_1 \symm F_3| \leq 2\alpha \sqrt{8|F_2| + 1} - 2 \alpha.
\]
\end{proof}

While we may not be able to reconstruct the set $F_2$, as it is not uniquely determined, we can reconstruct a uniquely determined neighbour $F_1$ of $F_2$. When $F_2$ is quite close to being uniquely determined, it must have a large intersection with $F_1$. Hence we know that at least a certain fraction of the points of $F_1$ must belong to $F_2$. The next theorem gives a bound for this fraction.

\begin{theorem}\label{metF1}
Let $F_2$ be a subset of $\mathbb{Z}^2$. Let $F_1$ be a uniquely determined neighbour of $F_2$. Put $\alpha = \alpha(F_2, F_1)$. Then
\[
\frac{|F_2 \cap F_1|}{|F_2|} \geq 1 - \frac{\sqrt{2} \alpha}{\sqrt{|F_2|}}.
\]
\end{theorem}

\begin{proof}
By Remark \ref{wortelgrens2} we have $|F_1 \symm F_2| \leq 2 \alpha \sqrt{2|F_2|}$.  Hence
\[
|F_1 \cap F_2| = |F_2| - \frac{1}{2} |F_1 \symm F_2| \geq |F_2| - \alpha \sqrt{2|F_2|}.
\]
Dividing by $|F_2|$ yields the theorem.
\end{proof}

Similarly, we can find a lower bound on the part of $F_2$ that must belong to any other reconstruction $F_3$.

\begin{theorem}\label{beta}
Let $F_2$ and $F_3$ be finite subsets of $\mathbb{Z}^2$ with the same line sums. Let $F_1$ be a uniquely determined neighbour of $F_2$ and $F_3$. Put $\alpha = \alpha(F_2, F_1)$. Then
\[
\frac{|F_2 \cap F_3|}{|F_2|} \geq 1 - \frac{2 \sqrt{2} \alpha}{\sqrt{|F_2|}}.
\]
\end{theorem}

\begin{proof}
By Remark \ref{wortelgrens2} we have $|F_1 \symm F_2| \leq 2 \alpha \sqrt{2|F_2|}$ and $|F_1 \symm F_3| \leq 2 \alpha \sqrt{2|F_2|}$. Hence
\[
|F_2 \symm F_3| \leq 4 \alpha \sqrt{2|F_2|}.
\]
So
\[
|F_2 \cap F_3| = |F_2| - \frac{1}{2} |F_2 \symm F_3| \geq |F_2| - 2 \alpha \sqrt{2|F_2|}.
\]
Dividing by $|F_2|$ yields the theorem.
\end{proof}

\begin{corollary}
If $F_2$ and $F_3$ are disjoint sets with the same line sums, then
\[
|F_2| \leq 8\alpha^2.
\]
\end{corollary}

\begin{proof}
If $F_2$ and $F_3$ are disjoint sets, then $\frac{|F_2 \cap F_3|}{|F_2|} = 0$, so by Theorem \ref{beta}
\[
0 \geq 1- \frac{2 \sqrt{2} \alpha}{\sqrt{|F_2|}},
\]
which we can rewrite as $|F_2| \leq 8\alpha^2$.
\end{proof}

Theorem \ref{metF1} shows that for given row and column sums that a set $F_2$ must satisfy, we can find a set of points $F_1$ such that a any possible set $F_2$ must contain a subset of $F_1$ with a certain size. However, it may happen that none of the individual points of $F_1$ must necessarily belong to such a set $F_2$. It is possible to determine from the line sums the intersection of all possible sets $F_2$, see e.g. \cite[Theorem 3.4]{anstee}. The following statement is a particular case of that theorem.

\begin{theorem}\label{notsure}
Let $F_2$ be a subset of $\mathbb{Z}^2$ with column sums $c_1^{(2)} \geq c_2^{(2)} \geq \ldots \geq c_n^{(2)}$. Let $F_1$ be a uniquely determined neighbour of $F_2$ with column sums $c_1^{(1)} \geq c_2^{(1)} \geq \ldots \geq c_n^{(1)}$. Suppose
\[
\sum_{j=1}^l c_j^{(1)} > \sum_{j=1}^l c_j^{(2)} \quad \text{for } l = 1, 2, \ldots, n-1.
\]
Then for all $(i,j) \in F_2$ there exists a set $F_3$ with the same row and column sums as $F_2$ such that $(i,j) \not\in F_3$.
\end{theorem}

We illustrate the theorems in this section by the following example.

\begin{example}\label{ex1}
Let $m$ and $n$ be positive integers. Let row sums $r_1$, $r_2$, \ldots, $r_n$ be given by $r_i = (n-i+1)m$ for $1 \leq i \leq n$. Let column sums $c_1$, $c_2$, \ldots, $c_{(n+1)m}$ be given by
\begin{itemize}
\item $c_j = n-1$ for $1 \leq j \leq m$,
\item $c_{lm+j} = n-l$ for $1 \leq l \leq n-1$, $1 \leq j \leq m$.
\item $c_j = 1$ for $nm+1 \leq j \leq (n+1)m$.
\end{itemize}
The uniquely determined set $F_1$ with row sums $r_1$, $r_2$, \ldots, $r_n$ has column sums $c_1'$, $c_2'$, \ldots, $c_{(n+1)m}'$ given by $c_{lm+j}' = n-l$ for $0 \leq l \leq n$, $1 \leq j \leq m$. For any set $F_2$ with row sums $r_1$, $r_2$, \ldots $r_n$ and column sum $c_1$, $c_2$, \ldots, $c_{(n+1)m}$ we have $\alpha = \alpha(F_1, F_2) = m$: the row sums of $F_1$ and $F_2$ are the same, while the column sums of the first $m$ and last $m$ columns differ by exactly 1.

Construct sets $F_2$ and $F_3$ as follows. In row $i$, $1 \leq i \leq n$, the elements of $F_2$ are in columns 1, 2, \ldots, $(n-i)m$ and in columns $(n-i+1)m + 1$, $(n-i+1)m +2$, \ldots, $(n-i+2)m$. In row $i$, $1 \leq i \leq n-1$, the elements of $F_3$ are in columns 1, 2, \ldots, $(n-i+1)m$. In row $n$ the elements of $F_3$ are in columns $nm+1$, $nm+2$, \ldots, $(n+1)m$. The sets $F_2$ and $F_3$ both have row sums $r_1$, $r_2$, \ldots $r_n$ and column sum $c_1$, $c_2$, \ldots, $c_{(n+1)m}$. We have $|F_2| = |F_3| = |F_1| = \frac{1}{2}mn(n+1)$.

\begin{figure}
\begin{center}
\includegraphics{plaatjesdoorsneden.1}
\end{center}
\caption{Example \ref{ex1} with $n=5$ and $m=3$. The set $F_2$ consists of the white and black-and-white points, while $F_3$ consists of the black and black-and-white points.}
\end{figure}

Theorem \ref{dezelfde} states that
\[
|F_2 \symm F_3| \leq 2m \sqrt{4mn(n+1)+1} -2m,
\]
while it actually holds that $|F_2 \symm F_3| = 2mn$.

Theorem \ref{metF1} states that
\[
\frac{|F_1 \cap F_2|}{|F_2|} \geq 1 - \frac{\sqrt{2}m}{\sqrt{\frac{1}{2}mn(n+1)}} \geq 1 - \frac{2\sqrt{m}}{n},
\]
while it actually holds that
\[
\frac{|F_1 \cap F_2|}{|F_2|} = \frac{\frac{1}{2}mn(n-1)}{\frac{1}{2}mn(n+1)} = \frac{n-1}{n+1} = 1 - \frac{2}{n+1}.
\]

Finally note that $F_2$ meets the conditions of Theorem \ref{notsure}, so none of the points of $F_2$ is contained in every set that has the same line sums as $F_2$.
\end{example}

\section{Sets with different line sums}\label{differentlinesums}

First consider two uniquely determined finite subsets $F_1$ and $F_1'$ of $\mathbb{Z}^2$. Let the row sums of $F_1$ be denoted by $r_1$, $r_2$, \ldots, $r_m$ and let the row sums of $F_1'$ be denoted by $r_1'$, $r_2'$, \ldots, $r_m'$. Without loss of generality, we may assume that $r_1 \geq r_2 \geq \ldots \geq r_m$.

Let $T$ be a staircase of which the elements are contained in the rows $i_1 < i_2 < \ldots < i_k$. Let $(i_t,j) \in F_1 \backslash F_1'$ and $(i_{t+1}, j) \in F_1' \backslash F_1$ be elements of $T$. By Lemma \ref{lemmauniek} we have $r_{i_t} > r_{i_{t+1}}$ and $r_{i_t}' < r_{i_{t+1}}'$. Row $i_1$ must contain an element of $F_1 \backslash F_1'$ of $T$, and row $i_k$ must contain an element of $F_1' \backslash F_1$ of $T$. Hence we can apply this for $t = 1, 2, \ldots, k-1$, and we find
\[
r_{i_1} > r_{i_2} > \ldots > r_{i_k},
\]
\[
r_{i_1}' < r_{i_2}' < \ldots < r_{i_k}'.
\]
Assume without loss of generality that there is at least one value of $t$ for which $r_{i_t}' - r_{i_t} \geq 0$. (Otherwise, reverse the role of $r_i'$ and $r_i$ in what follows.) Let \[u = \min\{r_{i_t}'-r_{i_t} :  r_{i_t}' - r_{i_t} \geq 0\}\] and let $s$ be such that $r_{i_s}' - r_{i_s} = u$. We distinguish two cases: $u =0$ and $u \geq 1$.

Case 1: suppose $u=0$. For $t \geq s$ we have $r_{i_t} \leq r_{i_s} - (t-s)$ and $r_{i_t}' \geq r_{i_s}' + (t-s)$, hence
\[
r_{i_t}' - r_{i_t} \geq r_{i_s}' - r_{i_s} + 2(t-s) = 2(t-s) \geq 0,
\]
so
\[
|r_{i_t}' - r_{i_t}| \geq 2(t-s).
\]
For $t < s$ we have $r_{i_t} \geq r_{i_s} + (s-t)$ and $r_{i_t}' \leq r_{i_s}' - (s-t)$, hence
\[
r_{i_t}' - r_{i_t} \leq r_{i_s}' - r_{i_s} - 2(s-t) = -2(s-t) < 0,
\]
so
\[
|r_{i_t}' - r_{i_t}| \geq 2(s-t).
\]
Now we have
\begin{eqnarray*}
\sum_{t=1}^k |r_{i_t}' - r_{i_t}| & \geq & \sum_{t=1}^{s-1} 2(s-t) + \sum_{t=s}^k 2(t-s) \\ & = & 2s^2 + (-2k -2)s + (k^2 +k) \\ & \geq & 2\left(\frac{k+1}{2}\right)^2 + (-2k-2) \frac{k+1}{2} + (k^2 + k) \\ & = & \tfrac{1}{2}k^2 - \tfrac{1}{2}.
\end{eqnarray*}

Case 2: suppose $u \geq 1$. Similarly to the first case, we have for $t \geq s$:
\[
|r_{i_t}' - r_{i_t}| \geq 2(t-s) + 1.
\]
If $s = 1$, there are no $t < s$ to consider. Assume $s \geq 2$. Then $r_{i_{s-1}}' - r_{i_{s-1}} < r_{i_s}' - r_{i_s} = u$, so by the minimality of $u$ we must have $r_{i_{s-1}}' - r_{i_{s-1}} \leq -1$. Similarly to above, we have
\[
|r_{i_t}' - r_{i_t}| \geq 2(s-t) - 1.
\]
Hence
\begin{eqnarray*}
\sum_{t=1}^k |r_{i_t}' - r_{i_t}| & \geq & \sum_{t=1}^{s-1} ( 2(s-t) -1 ) + \sum_{t=s}^k (2(t-s) +1) \\ & = & 2s^2 + (-2k-4)s + (k^2 +2k +2) \\ & \geq & 2\left(\frac{k+2}{2}\right)^2 + (-2k-4) \frac{k+2}{2} + (k^2 + 2k + 2) \\ & = & \tfrac{1}{2}k^2.
\end{eqnarray*}

In both cases we have $\sum_{t=1}^k |r_{i_t}' - r_{i_t}| \geq \tfrac{1}{2}k^2 - \tfrac{1}{2}$, and since the sum must be an integer, we have
\[
\sum_{t=1}^k |r_{i_t}' - r_{i_t}| \geq \lfloor \tfrac{1}{2}k^2 \rfloor.
\]
Hence the difference between the row sums of $F_1$ and $F_1'$ is at least $\lfloor \frac{1}{2}k^2 \rfloor$. Similarly, if $T$ is a staircase that contains elements in $k$ columns, the difference between the column sums of $F_1$ and $F_1'$ is at least $\lfloor \frac{1}{2} k^2 \rfloor$.

\begin{theorem}\label{uniek}
Let $F_1$ and $F_1'$ be uniquely determined finite subsets of $\mathbb{Z}^2$. Put $\alpha_1 = \alpha(F_1, F_1')$. Then
\[
|F_1 \symm F_1'| \leq 2\alpha_1 \sqrt{2\alpha_1 + 1} - \alpha_1.
\]
\end{theorem}

\begin{proof}
Consider all staircases in $F_1 \symm F_1'$, and let $T$ be one with the maximal number of elements. We distinguish two cases.
\begin{itemize}
\item Suppose $T$ has $2k+1$ elements for some $k \geq 0$. Then exactly $k+1$ rows and $k+1$ columns contain elements of $T$. By the argument above, we have
    \[
    2 \alpha_1 \geq \left\lfloor \tfrac{1}{2}(k+1)^2 \right\rfloor + \left\lfloor \tfrac{1}{2}(k+1)^2 \right\rfloor \geq (k+1)^2 - 1 = k^2 +2k.
    \]
    This implies $k +1 \leq \sqrt{2\alpha_1 + 1}$ and therefore $2k+1 \leq 2\sqrt{2\alpha_1 + 1} -1$.
\item Suppose $T$ has $2k$ elements for some $k \geq 1$. Then either $k$ rows and $k+1$ columns or $k+1$ rows and $k$ columns contain elements of $T$. By the argument above, we have
    \[
    2 \alpha_1 \geq \left\lfloor \tfrac{1}{2}(k+1)^2 \right\rfloor + \left\lfloor \tfrac{1}{2}k^2 \right\rfloor = \tfrac{1}{2}(k+1)^2 + \tfrac{1}{2}k^2 - \tfrac{1}{2} = k^2 + k.
    \]
    This implies $k+\frac{1}{2} \leq \sqrt{2\alpha_1 + \frac{1}{4}}$ and therefore $2k \leq 2 \sqrt{2\alpha_1+ \frac{1}{4}} - 1$.
\end{itemize}
All $\alpha_1$ staircases of $F_1 \symm F_2$ have at most as many elements as $T$, so in both cases we have
\[
|F_1 \symm F_1'| \leq 2 \alpha_1 \sqrt{2\alpha_1 + 1} - \alpha_1.
\]
\end{proof}

\begin{remark}
It is remarkable that the bound in Theorem \ref{uniek} does not depend on the sizes of $F_1$ and $F_1'$. Such a dependency cannot be avoided if one of the two sets is not uniquely determined, as in Lemma \ref{wortelgrens}. To show this, notice that in Example \ref{ex1} for fixed $\alpha = m$ the symmetric difference $|F_1 \symm F_2|$ becomes arbitrarily large when $n$ tends to infinity. Theorem \ref{uniek} shows that this cannot happen if both sets are uniquely determined.
\end{remark}

\begin{example}\label{ex2}
Let $n > 1$ be an integer. Define $r_i = n-i$ for $1 \leq i \leq n$ and $r_n' = n$. Let $F_1$ be the uniquely determined set with row and column sums $r_1$, $r_2$, \ldots, $r_n$. Let $F_1'$ be the uniquely determined set with row and column sums $r_1$, $r_2$, \ldots, $r_{n-1}$, $r_n'$. We have $\alpha_1 = \alpha(F_1, F_1') = n$. Consider row $i$, where $1 \leq i \leq n-1$. The elements of $F_1$ in this row are in columns 1, 2, \ldots, $n-i$, while the elements of $F_1'$ in this row are in columns 1, 2, \ldots, $n-i-1$ and $n$. In row $n$ there are $n$ elements of $F_1'$ and none of $F_1$. Hence
\[
|F_1 \symm F_1'| = 2(n-1) + n = 3n-2,
\]
while Theorem \ref{uniek} states that
\[
|F_1 \symm F_1'| \leq 2n \sqrt{2n + 1} - n.
\] \hfill $\square$

\begin{figure}
\begin{center}
\includegraphics{plaatjesdoorsneden.2}
\end{center}
\caption{Example \ref{ex2} with $n=7$. The set $F_1$ consists of the white and black-and-white points, while $F_1'$ consists of the black and black-and-white points.}
\end{figure}

\end{example}

Finally we derive a bound on the symmetric difference of two sets $F_2$ and $F_3$ with arbitrary line sums.

\begin{theorem}\label{verschillend}
Let $F_2$ and $F_3$ be finite subsets of $\mathbb{Z}^2$. Let $F_1$ be a uniquely determined neighbour of $F_2$, and let $F_1'$ be a uniquely determined neighbour of $F_3$. Put $\alpha_2 = \alpha(F_1,F_2)$, $\alpha_3 = \alpha(F_1', F_3)$ and $\alpha_1 = \alpha(F_1, F_1')$. Then
\[
|F_2 \symm F_3| \leq \alpha_2 \sqrt{8|F_2| + 1} - \alpha_2 + \alpha_3 \sqrt{8|F_3|+1} - \alpha_3 + 2\alpha_1 \sqrt{2\alpha_1 + 1} - \alpha_1.
\]
\end{theorem}

\begin{proof}
This is an immediate result of Lemma \ref{wortelgrens} and Theorem \ref{uniek}.
\end{proof}

\begin{example}\label{ex3}
Let $n$ be a positive integer. We construct sets $F_2$ and $F_3$ as follows.
\begin{itemize}
\item In row $i$, where $1 \leq i \leq n$, the elements of $F_2$ are in columns 1, 2, \ldots, $2(n-i)$ as well as columns $2(n-i)+2$ and $2(n-i)+3$.
\item In row $n+1$, there is a single element of $F_2$ in column 1.
\item In row $i$, where $1 \leq i \leq n$, the elements of $F_3$ are in columns 1, 2, \ldots, $2(n-i)+1$ as well as column $2(n-i)+4$.
\item In row $n+1$ there are no elements of $F_3$.
\end{itemize}
The row sums of $F_2$ are given by
\[
r_i^{(2)} = \left\{ \begin{array}{cl} 2(n-i+1) & \text{if } 1 \leq i \leq n, \\ 1 & \text{if } i=n+1. \end{array} \right.
\]
The column sums of $F_2$ are given by
\[
c_j^{(2)} = \left\{ \begin{array}{cl} n - \lfloor \frac{j-1}{2} \rfloor & \text{if } 1 \leq j \leq 2n, \\ 1 & \text{if } j=2n+1, \\ 0 & \text{if } j=2n+2. \end{array} \right.
\]
The row sums of $F_3$ are given by
\[
r_i^{(3)} = 2(n-i+1), \quad 1 \leq i \leq n+1.
\]
The column sums of $F_3$ are given by
\[
c_j^{(3)} = \left\{ \begin{array}{cl} n & \text{if } j=1, \\ n-1 & \text{if } j=2, \\ n - \lfloor \frac{j-1}{2} \rfloor & \text{if } 3 \leq j \leq 2n, \\ 0 & \text{if } j=2n+1, \\ 1 & \text{if } j=2n+2. \end{array} \right.
\]

\begin{figure}
\begin{center}
\includegraphics{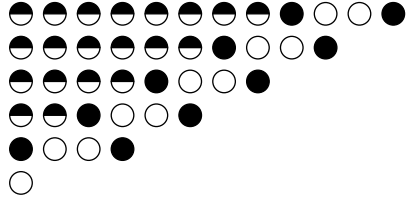}
\end{center}
\caption{Example \ref{ex3} with $n=5$. The set $F_2$ consists of the white and black-and-white points, while $F_3$ consists of the black and black-and-white points.}
\end{figure}

Let $F_1$ be the uniquely determined set with the same row sums as $F_2$ and non-increasing column sums. Let $F_1'$ be the uniquely determined set with the same row sums as $F_3$ and non-increasing column sums. We have
\[
\alpha_2 = \alpha(F_2, F_1) = 1,
\quad
\alpha_3 = \alpha(F_3, F_1') = 1,
\quad
\alpha_1 = \alpha(F_1, F_1') = 1.
\]
Furthermore, $|F_2| = n(n+1)+1$ and $|F_3| = n(n+1)$.

Theorem \ref{verschillend} states that
\[
|F_2 \symm F_3| \leq  \sqrt{8n(n+1) + 9} + \sqrt{8n(n+1)+1} + 2 \sqrt{3} - 3 \approx 4\sqrt{2} n,
\]
while actually
\[
|F_2 \symm F_3| = 4n+1.
\]
\end{example}

\section{Concluding remarks}

We have proved an upper bounds on the difference between two images with the same row and column sums, as well as on the difference between two images with different row and column sums. The bounds heavily depend on the parameter $\alpha$, which indicates how close an image is to being uniquely determined. If a set of given line sums ``almost uniquely determines'' the image (i.e. $\alpha$ is very small) it may still happen that no points belong to all possible images with those line sums. However, using the results from this paper we can find a set of points of which a subset of certain size is guaranteed to belong to the image.

There is still a gap between the examples we have found and the bounds we have proven. It appears that all bounds can be improved by a factor $\sqrt{\alpha}$. For this it would suffice to improve both Lemma \ref{wortelgrens} and Theorem \ref{uniek} by a factor $\sqrt{\alpha}$, but so far we did not manage to improve either of those.

The results of this paper can be applied to projections in more than two directions as well: simply pick two directions and forget about the others. One would expect this to give bad results, but that is actually not always the case. It is possible to construct examples with projections in more than two directions where the bound using only two of the directions is still only a factor $\sqrt{\alpha}$ off. However, in many cases it should be (somehow) possible to use the projections in all directions to get better results. Our future research will be concerned with finding such a method.

\end{document}